\documentclass[a4paper,11pt,twoside,english]{article}
\usepackage{babel}
\usepackage{latexsym,amsfonts}
\usepackage{amsthm}
\usepackage{amsmath}
\usepackage{mathtools}
\usepackage{paralist}
\usepackage{ifthen}
\usepackage{moredefs}
\usepackage{todonotes}

\usepackage{lastpage}
\usepackage{fancyhdr}
\fancypagestyle{plain}{

  \fancyhead[L,C]{}
  \fancyhead[R]{\today}
  \fancyfoot[R,L]{}
  \fancyfoot[C]{\thepage \ of \pageref{LastPage}}
  }

\title{WORTH property, Garc\'{i}a-Falset coefficient and Opial property of infinite sums}
\author{Jan-David Hardtke}
\date{}

\DeclareMathOperator{\lin}{span}

\DeclareMathOperator{\diam}{diam}

\providecommand{\shortversion}[1]{#1}
\providecommand{\longversion}[2]{#2 {\rm(#1 in short)}}

\providecommand{\acronym}[2]{
  \newboolean{acronym#1}
  \setboolean{acronym#1}{true}
  \expandafter\providecommand\expandafter{\csname acronym#1\endcsname}{#2}
  }

\providecommand{\ac}[1]{\ifthenelse{\boolean{acronym#1}}
{\longversion{#1}{\csname acronym#1\endcsname}\Global\ToggleBoolean{acronym#1}}
{\shortversion{#1}}}

\providecommand{\ifif}{iff }
\providecommand{\sm}{\setminus}
\providecommand{\ssq}{\subseteq}

\providecommand{\N}{\ensuremath{\mathbb{N}}}

\providecommand{\R}{\ensuremath{\mathbb{R}}}

\providecommand{\eps}{\ensuremath{\varepsilon}}

\providecommand{\keywords}[1]{
  {\let\thefootnote=\relax
  \footnote{{\em Keywords}: #1}}
  \addtocounter{footnote}{-1}
  }

\providecommand{\AMS}[1]{
  {\let\thefootnote=\relax
  \footnote{{\em AMS Subject Classification} (2010): #1}}
  \addtocounter{footnote}{-1}
  }

\providecommand{\address}{
  {\sc \noindent Department of Mathematics \\
  Freie Universit\"at Berlin \\  
  Arnimallee 6, 14195 Berlin \\
  Germany \\}
  }

\DeclarePairedDelimiter{\set}{\lbrace}{\rbrace}
\DeclarePairedDelimiter{\paren}{\lparen}{\rparen}

\DeclarePairedDelimiter{\abs}{\lvert}{\rvert}
\DeclarePairedDelimiter{\norm}{\lVert}{\rVert}

\theoremstyle{definition}
\newtheorem{definition}{Definition}[section]
\newtheorem*{definition*}{Definition}
\newtheorem{example}[definition]{Example}
\newtheorem*{example*}{Example}

\newtheorem*{remark*}{Remark}

\newtheorem*{problem*}{Problem}

\newtheorem*{question*}{Question}

\newtheorem*{conjecture*}{Conjecture}
\theoremstyle{remark}

\newtheorem*{claim*}{Claim}

\newtheorem*{fact*}{Fact}
\theoremstyle{plain}
\newtheorem{lemma}[definition]{Lemma}
\newtheorem*{lemma*}{Lemma}
\newtheorem{proposition}[definition]{Proposition}
\newtheorem*{proposition*}{Proposition}
\newtheorem{theorem}[definition]{Theorem}
\newtheorem*{theorem*}{Theorem}
\newtheorem{corollary}[definition]{Corollary}
\newtheorem*{corolary*}{Corollary}

\newenvironment{Proof}[1][\proofname]{\begin{proof}[#1] \setlength{\parindent}{0pt}}{\end{proof}}
\newenvironment{Abstract}{\centering\begin{minipage}{0.8\textwidth} \noindent \small {\sc Abstract.}}{\end{minipage}\par}

\usepackage{color}
\definecolor{darkgreen}{rgb}{0,0.5,0}

\numberwithin{equation}{section}
\newtagform{colored}[\color{blue}]{\color{blue}(}{\color{blue})}
\usetagform{colored}

\usepackage[colorlinks,linkcolor=blue,citecolor=red,urlcolor=darkgreen]{hyperref}
\providecommand{\email}{{\it E-mail address:} \href{mailto:hardtke@math.fu-berlin.de}{\tt hardtke@math.fu-berlin.de}}

\usepackage{amsrefs}

\begin{document}

\maketitle

\begin{Abstract}
We prove some results concerning the WORTH pro\-perty and the Garc\'{i}a-Falset coefficient of absolute sums of infinitely 
many Banach spaces. The Opial property/uniform Opial property of infinite $\ell^p$-sums is also studied and some properties 
analogous to the Opial property/uniform Opial property for Lebesgue-Bochner spaces $L^p(\mu,X)$ are discussed.
\end{Abstract}
\AMS{46B20; 46E40}

\section{Introduction}\label{sec:intro}
For a real Banach space $X$, denote by $X^*$ its dual space, by $B_X$ its closed unit ball and by $S_X$ its unit sphere.\par
We begin by recalling the important notion of fixed point property: $X$ is said to have the fixed point property 
(resp. weak fixed point property) if for every closed and bounded (resp. weakly compact) convex subset $C\ssq X$, 
every nonexpansive mapping $F:C \rightarrow C$ has a fixed point (where $F$ is called nonexpansive if 
$\norm{F(x)-F(y)}\leq\norm{x-y}$ for all $x,y\in C$, in other words, if $F$ is $1$-Lipschitz continuous).\par
A bounded closed convex subset $C\ssq X$ is said to have normal structure if for each subset $B\ssq C$ which contains
at least two elements there exists a point $x\in B$ such that
\begin{equation*}
\sup_{y\in B}\norm{x-y}<\diam{B}.
\end{equation*}
It is well known that if $C$ is weakly compact and has normal structure, then every nonexpansive mapping 
$F:C \rightarrow C$ has a fixed point (see for example \cite{goebel2}*{Theorem 2.1}).\par
The space $X$ is said to have the Opial property provided that 
\begin{equation*}
\limsup_{n\to \infty}\norm{x_n}<\limsup_{n\to \infty}\norm{x_n-x}
\end{equation*}
holds for every weakly nullsequence $(x_n)_{n\in \N}$ in $X$ and every $x\in X\setminus\set*{0}$ (one could
as well use $\liminf$ instead of $\limsup$ or assume from the beginning that both limits exist). This property 
was first considered by Opial in \cite{opial} (starting from the Hilbert spaces as canonical example) to provide 
a result on iterative approximations of fixed points of nonexpansive mappings. It is shown in \cite{opial} that
the spaces $\ell^p$ for $1\leq p<\infty$ enjoy the Opial property, whereas $L^p[0,1]$ for $1<p<\infty, p\neq 2$
fails to have it. Note further that every Banach space with the Schur property (i.\,e. weak and norm convergence of 
sequences coincide) trivially has the Opial property. Also, $X$ is said to have the nonstrict Opial property if it fulfils 
the definition of the Opial property with ``$\leq$'' instead of ``$<$'' (\cite{sims2}, in \cite{garcia-falset1} it is called 
weak Opial property). It is known that every weakly compact convex set in a Banach space with the Opial property has normal 
structure (see for instance \cite{prus}*{Theorem 5.4}).\par
Prus introduced the notion of uniform Opial property in \cite{prus2}: a Banach space $X$ has the uniform Opial property if
for every $c>0$ there is some $r>0$ such that
\begin{equation*}
1+r\leq\liminf_{n\to \infty}\norm{x_n-x}
\end{equation*}
holds for every $x\in X$ with $\norm{x}\geq c$ and every weakly nullsequence $(x_n)_{n\in \N}$ in $X$ with 
$\liminf\norm{x_n}\geq 1$. In \cite{prus2} it was proved that a Banach space is reflexive and has the uniform
Opial property if and only if it has the so called property $(L)$ (see \cite{prus2} for the definition), and that 
$X$ has the fixed point property whenever $X^*$ enjoys said property $(L)$.\par
A modulus corresponding to the uniform Opial property was defined in \cite{lin}:
\begin{equation*}
r_X(c):=\inf\set*{\liminf_{n\to \infty}\norm{x_n-x}-1} \ \ \forall c>0,
\end{equation*}
where the infimum is taken over all $x\in X$ with $\norm{x}\geq c$ and all weakly nullsequences $(x_n)_{n\in \N}$ in $X$ with 
$\liminf\norm{x_n}\geq 1$ (if $X$ has the Schur property, we agree to set $r_X(c):=1$ for all $c>0$). Then $X$ has the uniform 
Opial property \ifif $r_X(c)>0$ for every $c>0$.\par
In this paper, we will mostly use the following equivalent formulation of the uniform Opial property (\cite{khamsi}*{Definition 3.1}):
$X$ has the uniform Opial property \ifif for every $\eps>0$ and every $R>0$ there is some $\eta>0$ such that
\begin{equation*}
\eta+\liminf_{n\to \infty}\norm{x_n}\leq\liminf_{n\to \infty}\norm{x_n-x}
\end{equation*}
holds for all $x\in X$ with $\norm{x}\geq\eps$ and every weakly nullsequence $(x_n)_{n\in \N}$ in $X$ with $\limsup\norm{x_n}\leq R$.\par
We can also associate a modulus to this formulation in the following way:
\begin{equation*}
\eta_X(\eps,R):=\inf\set*{\liminf_{n\to \infty}\norm{x_n-x}-\liminf_{n\to \infty}\norm{x_n}} \ \ \forall \eps,R>0,
\end{equation*}
where the infimum is taken over all $x\in X$ with $\norm{x}\geq\eps$ and all weakly nullsequences $(x_n)_{n\in \N}$ in $X$ 
with $\limsup\norm{x_n}\leq R$. So $X$ has the uniform Opial property \ifif $\eta_X(\eps,R)>0$ for all $\eps,R>0$. Actually,
it is enough that for every $\eps>0$ there exists some $R>2$ with $\eta_X(\eps,R)>0$. More precisely, we have the following connection 
between the two moduli $r_X$ and $\eta_X$.
\begin{lemma}\label{lemma:opialmodul}
Let $X$ be a Banach space which does not have the Schur property.
\begin{enumerate}[\upshape(i)]
\item For every $c>0$ and every $R>2$ we have
\begin{equation*}
\min\set*{\eta_X(c,R),\frac{R}{2}-1}\leq r_X(c).
\end{equation*}
\item For all $\eps,R>0$ with $r_X(\eps/R)>0$ we have
\begin{equation*}
\frac{\eps r_X(\eps/R)}{2+r_X(\eps/R)}=\max_{\beta\in [0,\eps/2]}\min\set*{\beta r_X(\frac{\eps}{R}),\eps-2\beta}\leq \eta_X(\eps,R).
\end{equation*}
\end{enumerate}
\end{lemma}

\begin{Proof}
(i) Let $c>0$ and $R>2$. Put $\tau:=\min\set*{\eta_X(c,R),\frac{R}{2}-1}$. Let $(x_n)_{n\in \N}$ be any weakly nullsequence 
in $X$ with $\liminf\norm{x_n}\geq 1$ and let $x\in X$ with $\norm{x}\geq c$. By passing to a subsequence, we may assume that 
$\lim_{n\to \infty}\norm{x_n-x}$ and $s:=\lim_{n\to \infty}\norm{x_n}$ exist. If $s\leq R$ then $1+\tau\leq s+\eta_X(c,R)\leq
\lim_{n\to \infty}\norm{x_n-x}$.\par
If $s>R$ and $\norm{x}>R/2$, then $\lim_{n\to \infty}\norm{x_n-x}\geq\norm{x}>R/2\geq1+\tau$ by the weak lower semicontinuity 
of the norm. Finally, if $s>R$ and $\norm{x}\leq R/2$, then $\lim_{n\to \infty}\norm{x_n-x}\geq s-\norm{x}>R/2\geq1+\tau$.\par
(ii) The first equality is easily verified. Now chose any $\beta\in (0,\eps/2)$ and put $\nu:=\min\set*{\beta r_X(\frac{\eps}{R}),\eps-2\beta}$.
Let $(x_n)_{n\in \N}$ be a weakly nullsequence in $X$ with $\limsup\norm{x_n}\leq R$ and let $x\in X$ with $\norm{x}\geq\eps$. 
Again we may assume that $\lim_{n\to \infty}\norm{x_n-x}$ and $s:=\lim_{n\to \infty}\norm{x_n}$ exist. By the definition of $r_X$ we get 
$s(1+r_X(\eps/R))\leq\lim_{n\to \infty}\norm{x_n-x}$, which implies $s+\nu\leq\lim_{n\to \infty}\norm{x_n-x}$ if $s>\beta$. But if $s\leq\beta$
then $\lim_{n\to \infty}\norm{x_n-x}\geq\norm{x}-s\geq\eps-\beta\geq\nu+\beta\geq\nu+s$ and the proof is finished.
\end{Proof}

In \cite{garcia-falset1} J. Garc\'{i}a-Falset introduced the following coefficient of a Banach space $X$:
\begin{equation*}
R(X):=\sup\set*{\liminf_{n\to \infty}\norm{x_n+x}: x\in B_X, (x_n)_{n\in \N}\in \operatorname{WN}(B_X)},
\end{equation*}
where we denote by $\operatorname{WN}(B_X)$ the set of all weakly nullsequences in $B_X$. Obviously, $1\leq R(X)\leq 2$ 
and $R(X)=1$ if $X$ has the Schur property (in particular if $X$ is finite-dimensional or $X=\ell^1$). One has $R(c_0)=1$
and $R(\ell^p)=2^{1/p}$ for $1<p<\infty$ (see \cite{garcia-falset1}*{Corollary 3.2}). In \cite{garcia-falset2}*{Theorem 3}
it was proved that the condition $R(X)<2$ implies that $X$ has the weak fixed point property. The reflexive spaces with $R(X)<2$ 
are precisely the so called weakly nearly uniformly smooth spaces (\cite{garcia-falset1}*{Corollary 4.4}), which were introduced 
in \cite{kutzarova} and include in particular all uniformly smooth spaces.\par
We will denote by $\delta_X$ the modulus of convexity of $X$, i.\,e. for $0<\eps\leq 2$
\begin{equation*}
\delta_X(\eps):=\inf\set*{1-\norm*{\frac{x+y}{2}}:x,y\in B_X \ \text{with} \ \norm{x-y}\geq\eps}.
\end{equation*}
$X$ is uniformly rotund \ifif $\delta_X(\eps)>0$ for each $0<\eps\leq 2$. It is well-known that all spaces $L^p(\mu)$ for
any measure $\mu$ and any $1<p<\infty$ (in particular the spaces $\ell^p(I)$ for any index set $I$) are uniformly rotund.\par
In \cite{sims1} Sims introduced the notion of WORTH (weak orthogonality) property: $X$ is said to have the WORTH property provided that
for all weakly nullsequences $(x_n)_{n\in \N}$ in $X$ and every $x\in X$ one has $\norm{x_n+x}-\norm{x_n-x}\to 0$.\par
Again spaces with the Schur property obviously enjoy the WORTH property. Hilbert spaces are easily seen to have the WORTH property as well.
Also, the class of spaces with the WORTH property includes all so called weakly orthogonal Banach lattices (a notion introduced earlier by Borwein 
and Sims in \cite{borwein}), which in turn includes in particular all spaces $\ell^p(I)$ for $1\leq p<\infty$ and $c_0(I)$. However, the spaces
$L^p[0,1]$ with $1\leq p\leq\infty, p\neq 2$ do not have the WORTH property (see the remark at the end of \cite{sims2}). In \cite{sims1} it was 
proved that the WORTH property implies the nonstrict Opial property, and in \cite{sims2} it was shown that a space with the WORTH property which 
is $\eps$-inquadrate in every direction for some $0<\eps<2$ (see \cite{sims2} for the definition) has the weak fixed point property (even more, 
every weakly compact convex subset of such a space has normal structure). By \cite{garcia-falset1}*{Proposition 3.6}, a uniformly 
non-square\footnote{Recall that $X$ is said to be uniformly non-square if there is some $\delta>0$ such that whenever $x,y\in B_X$ one 
has $\norm{x+y}<2(1-\delta)$ or $\norm{x-y}<2(1-\delta)$.} Banach space $X$ with the WORTH property satisfies $R(X)<2$.\par
The degree $w(X)$ of WORTHness of $X$ was also introduced in \cite{sims2} as the supremum of all $r\geq 0$ such that
\begin{equation*}
r\liminf_{n\to \infty}\norm{x_n+x}\leq\liminf_{n\to \infty}\norm{x_n-x} 
\end{equation*}
holds for all $x\in X$ and all weakly nullsequences $(x_n)_{n\in \N}$ in $X$. Then $1/3\leq w(X)\leq 1$ and $X$ has the WORTH property 
if and only if $w(X)=1$.\par
In this paper, we will study the WORTH property and the Garc\'{i}a-Falset coefficient for infinite absolute sums, and the different Opial 
properties specifically for infinite $\ell^p$-sums of Banach spaces (for normal structure in (finite and infinite) direct sums of Banach 
spaces see \cite{dominguez-benavides} and references therein). The next section contains the necessary preliminaries on absolute sums.

\section{Preliminaries on absolute sums}\label{sec:sums}
Throughout this paper, if not otherwise stated, $I$ denotes a (mostly infinite) index set and $E$ a 
subspace of the space of all real-valued functions on $I$ which contains all functions with finite 
support and is endowed with an absolute, normalised norm $\norm{\cdot}_E$. The latter means that 
$\norm{\cdot}_E$ is a complete norm on $E$ such that the following conditions are satisfied:
\begin{enumerate}[(i)]
\item If $(a_i)_{i\in I}\in E$ and $(b_i)_{i\in I}\in \R^I$ such that $\abs{a_i}=\abs{b_i}$ for all
$i\in I$ then $(b_i)_{i\in I}\in E$ and $\norm{(b_i)_{i\in I}}_E=\norm{(a_i)_{i\in I}}_E$.
\item $\norm{e_i}_E=1$ for all $i\in I$, where $e_i=(e_{ij})_{j\in I}$ with $e_{ij}=0$ for $j\neq i$ 
and $e_{ii}=1$.
\end{enumerate}
It is important to note that such norms are automatically monotone, i.\,e. we actually have
\begin{align*}
&(a_i)_{i\in I}\in E, (b_i)_{i\in I}\in \R^I \ \ \text{with} \ \ \abs{b_i}\leq\abs{a_i} \ \forall i\in I \\ 
&\Rightarrow \ (b_i)_{i\in I}\in E \ \ \text{and} \ \ \norm{(b_i)_{i\in I}}_E\leq\norm{(a_i)_{i\in I}}_E.
\end{align*}
For a proof see for instance \cite{lee}*{Remark 2.1}. Standard examples of spaces with absolute, normalised 
norm are of course the spaces $\ell^p(I)$ (with $1\leq p\leq\infty$) and $c_0(I)$.\par
If we put 
\begin{equation*}
E^\prime:=\set*{(a_i)_{i\in I}\in \R^I: \norm{(a_i)_{i\in I}}_{E^{\prime}}:=\sup_{(b_i)_{i\in I}\in B_E}\sum_{i\in I}\abs*{a_ib_i}<\infty},
\end{equation*}
then $(E^{\prime},\norm{\cdot}_{E^{\prime}})$ is again a space with absolute, normalised norm and the map $T: E^\prime \rightarrow E^*$ defined by
\begin{equation*}
T((a_i)_{i\in I})((b_i)_{i\in I}):=\sum_{i\in I}a_ib_i \ \ \forall (a_i)_{i\in I}\in E^\prime, \forall (b_i)_{i\in I}\in E
\end{equation*}
is an isometric embedding. $T$ is onto if $\lin\set*{e_i:i\in I}$ is dense in $E$, so in this case $E^*=E^{\prime}$\par
Now given a family $(X_i)_{i\in I}$ of Banach spaces, the absolute sum of $(X_i)_{i\in I}$ with respect to $E$ is defined 
as the space
\begin{equation*}
\Bigl[\bigoplus_{i\in I}X_i\Bigr]_E:=\set*{(x_i)_{i\in I}\in \prod_{i\in I}X_i: (\norm*{x_i})_{i\in I}\in E}
\end{equation*}
endowed with the norm $\norm{(x_i)_{i\in I}}_E:=\norm{(\norm{x_i})_{i\in I}}_E$. It is not hard to see that 
this sum is indeed a Banach space. For $E=\ell^p(I)$ one obtains the usual $p$-sum.\par
As regards the dual space of an absolute sum, the map
\begin{align*}
&S: \Bigl[\bigoplus_{i\in I}X_i^*\Bigr]_{E^\prime} \rightarrow \Bigl[\bigoplus_{i\in I}X_i\Bigr]_E^* \\
&S((x_i^*)_{i\in I})((x_i)_{i\in I}):=\sum_{i\in I}x_i^*(x_i)
\end{align*}
is an isometric embedding and it is onto if $\lin\set*{e_i:i\in I}$ is dense in $E$.\par

\section{Results on absolute sums}\label{sec:results}

\subsection{WORTH property of absolute sums}\label{sub:worth}
By \cite{kato}*{Theorem 4.7}, $w(X\oplus_E Y)=\min\set*{w(X),w(Y)}$ holds for all Banach spaces $X$ and $Y$ and every absolute, normalised norm 
$\norm{\cdot}_E$ on $\R^2$ (actually, the notion of $\psi$-direct sums is used in \cite{kato}, but it is an equivalent formulation (see section 2 
in \cite{kato})). In particular, $X\oplus_E Y$ has the WORTH property if and only if $X$ and $Y$ have the WORTH property (for this, see also 
\cite{kato}*{Theorem 4.2}). It is possible to generalise \cite{kato}*{Theorem 4.7} to sums of arbitrarily many Banach spaces under a mild condition on $E$.
\begin{proposition}\label{prop:worth}
If $\lin\set*{e_i:i\in I}$ is dense in $E$ and $(X_i)_{i\in I}$ is any family of Banach spaces then
\begin{equation*}
w\paren*{\Bigl[\bigoplus_{i\in I}X_i\Bigr]_E}=\inf\set*{w(X_i):i\in I}.
\end{equation*}
In particular, $\bigl[\bigoplus_{i\in I}X_i\bigr]_E$ has the WORTH property if and only if $X_i$ has the WORTH property for every $i\in I$.
\end{proposition}

\begin{Proof}
Let us write $X=\bigl[\bigoplus_{i\in I}X_i\bigr]_E$ and $s=\inf\set*{w(X_i):i\in I}$. We clearly have $w(X)\leq s$. Now let $x_n=(x_{n,i})_{i\in I}\in X$ 
for every $n\in \N$ such that $(x_n)_{n\in \N}$ converges weakly to zero and let $x=(x_i)_{i\in I}\in X$. Without loss of generality, we may assume that
the limits $a:=\lim_{n\to \infty}\norm{x_n+x}_E$ and $b:=\lim_{n\to \infty}\norm{x_n-x}_E$ exist.\par
Since $\lin\set*{e_i:i\in I}$ is dense in $E$, it is not hard to see that we actually have $(\norm{x_i})_{i\in I}=\sum_{i\in I}\norm{x_i}e_i$. So if $\eps>0$ 
is given, we find a finite set $J\ssq I$ such that
\begin{equation}\label{eq:1}
\norm*{\sum_{i\in I\sm J}\norm{x_i}e_i}_E=\norm*{(\norm{x_i})_{i\in I}-\sum_{i\in J}\norm{x_i}e_i}_E\leq\eps.
\end{equation}
By passing to an appropriate subsequence we may assume that the limits $a_i:=\lim_{n\to \infty}\norm{x_{n,i}+x_i}$ and $b_i:=\lim_{n\to \infty}\norm{x_{n,i}-x_i}$ 
exist for each $i\in J$. Since $(x_{n,i})_{n\in \N}$ is weakly convergent to zero in $X_i$ for every $i\in I$ it follows that $sa_i\leq b_i\leq s^{-1}a_i$
and consequently
\begin{equation}\label{eq:2}
\abs{a_i-b_i}\leq\frac{1-s}{s}b_i \ \ \forall i\in J.
\end{equation}
For every $n\in \N$ we have, because of \eqref{eq:1},
\begin{align*}
&\abs{\norm{x_n+x}_E-\norm{x_n-x}_E}\leq\norm{(\norm{x_{n,i}+x_i}-\norm{x_{n,i}-x_i})_{i\in I}}_E \\
&\leq\norm*{\sum_{i\in J}(\norm{x_{n,i}+x_i}-\norm{x_{n,i}-x_i})e_i}_E+2\norm*{\sum_{i\in I\sm J}\norm{x_i}e_i}_E \\
&\leq\norm*{\sum_{i\in J}(\norm{x_{n,i}+x_i}-\norm{x_{n,i}-x_i})e_i}_E+2\eps.
\end{align*}
So for $n\to \infty$ we obtain
\begin{equation*}
\abs{a-b}\leq\norm*{\sum_{i\in J}(a_i-b_i)e_i}_E+2\eps.
\end{equation*}
Taking \eqref{eq:2} into account we arrive at
\begin{equation*}
\abs{a-b}\leq\frac{1-s}{s}\norm*{\sum_{i\in J}b_ie_i}_E+2\eps.
\end{equation*}
But $\norm{\sum_{i\in J}\norm{x_{n,i}-x_i}e_i}_E\leq\norm{x_n-x}_E$ for each $n$, thus $\norm{\sum_{i\in J}b_ie_i}_E\leq b$ and hence
\begin{equation*}
\abs{a-b}\leq\frac{1-s}{s}b+2\eps.
\end{equation*}
Letting $\eps\to 0$ leaves us with $\abs{a-b}\leq(1-s)b/s$ which implies $sa\leq b$ and we are done.
\end{Proof}

\subsection{Garc\'{i}a-Falset coefficient of absolute sums}\label{sub:garciafalset}
In \cite{dhompongsa}*{Theorem 7} it was proved that $R((X_1\oplus X_2\oplus\dots\oplus X_n)_E)<2$ whenever $R(X_i)<2$ for $i=1,\dots,n$ 
and $\norm{\cdot}_E$ is any strictly convex, absolute, normalised norm on $\R^n$. For absolute sums of two Banach spaces a stronger result 
was obtained in \cite{kato}*{Theorem 3.6}: $R(X\oplus_E Y)<2$ provided that $R(X), R(Y)<2$ and $\norm{\cdot}_E$ is any absolute, normalised
norm on $\R^2$ with $\norm{\cdot}_E\neq\norm{\cdot}_1$. For infinite sums we have the following theorem (for $J\ssq I$ we denote by 
$\bigl[\bigoplus_{i\in J}X_i\bigr]_E$ the sum of the family whose $i$-th member is $X_i$ for $i\in J$ and $\set*{0}$ for $i\in I\sm J$).
\begin{theorem}\label{thm:garcia-falset}
If $I$ is an infinite index set, $E$ a subspace of $\R^I$ with absolute, normalised norm such that $\lin\set*{e_i:i\in I}$ is dense in $E$
and $(X_i)_{i\in I}$ is a family of Banach spaces with
\begin{equation}\label{eq:3}
\alpha:=\sup\set*{R\paren*{\Bigl[\bigoplus_{i\in J}X_i\Bigr]_E}:J\ssq I \ \mathrm{finite}}<2
\end{equation}
and $\delta_E((1-\alpha/2)^2)>0$, then $R\paren*{\bigl[\bigoplus_{i\in I}X_i\bigr]_E}<2$.
\end{theorem}

\begin{Proof}
Let us write $X=\bigl[\bigoplus_{i\in I}X_i\bigr]_E$ for short. It is well known that $\delta_E$ is continuous on $(0,2)$ (see for example \cite{goebel}*{Lemma 5.1}),
so we can find $0<\tau<(1-\alpha/2)^2$ with $\delta_E(\tau)>0$. Let $\gamma:=\sqrt{\tau}$ and choose $0<\eta<\min\set*{\delta_E(\tau),1/2-\gamma}$.\par
Suppose that $R(X)=2$. Then there would be a weakly null sequence $(x_n)_{n\in \N}=((x_{n,i})_{i\in I})_{n\in \N}$ in $B_X$ and an element $x=(x_i)_{i\in I}\in B_X$ 
such that $\lim_{n\to \infty}\norm{x_n+x}_E>2-\eta$. We may assume $\norm{x_n+x}_E>2-2\eta$ for all $n\in \N$. Since $\norm{(\norm{x_{n,i}}+\norm{x_i})_{i\in I}}_E
\geq\norm{x_n+x}_E$ and $\eta<\delta_E(\tau)$ it follows that
\begin{equation}\label{eq:4}
\norm{(\norm{x_{n,i}}-\norm{x_i})_{i\in I}}_E<\tau \ \ \forall n\in \N.
\end{equation}
Similarly, 
\begin{equation*}
4(1-\eta)<2\norm{x_n+x}_E\leq\norm{(\norm{x_{n,i}}+\norm{x_i}+\norm{x_{n,i}+x_i})_{i\in I}}_E\leq 4
\end{equation*}
and hence
\begin{equation}\label{eq:5}
\norm{(\norm{x_{n,i}}+\norm{x_i}-\norm{x_{n,i}+x_i})_{i\in I}}_E<2\tau \ \ \forall n\in \N.
\end{equation}
We further have $\norm{x}_E\geq\norm{x_n+x}_E-1>1-2\eta>2\gamma$. Since $(\norm{x_i})_{i\in I}=\sum_{i\in I}\norm{x_i}e_i$ we can find a finite set $J\ssq I$
such that
\begin{equation}\label{eq:6}
\norm*{\sum_{i\in J}\norm{x_i}e_i}_E>2\gamma.
\end{equation}
Put $y:=(x_i)_{i\in J}, y_n=(x_{n,i})_{i\in J}\in \bigl[\bigoplus_{i\in J}X_i\bigr]_E$ as well as $a:=\sum_{i\in J}\norm{x_i}e_i$ and $a_n:=\sum_{i\in J}\norm{x_{n,i}}e_i$.
By \eqref{eq:4} we have
\begin{equation}\label{eq:7}
\norm{a_n-a}_E\leq\norm{(\norm{x_{n,i}}-\norm{x_i})_{i\in I}}_E<\tau \ \ \forall n\in \N,
\end{equation}
which implies in particular $\abs{\norm{y}_E-\norm{y_n}_E}=\abs{\norm{a}_E-\norm{a_n}_E}<\tau$, hence $\norm{y_n}_E>\norm{y}_E-\tau>2\gamma-\tau>0$, by \eqref{eq:6}.\par
Furthermore, for every $n\in \N$,
\begin{equation*}
\abs{\norm{a_n+a}_E-\norm{y_n+y}_E}\leq\norm*{\sum_{i\in J}(\norm{x_{n,i}}+\norm{x_i}-\norm{x_{n,i}+x_i})e_i}_E,
\end{equation*}
so because of \eqref{eq:5} it follows that
\begin{equation*}
\abs{\norm{a_n+a}_E-\norm{y_n+y}_E}<2\tau \ \ \forall n\in \N.
\end{equation*}
Also, by \eqref{eq:7}, we have $\abs{\norm{a_n+a}_E-2\norm{y}_E}=\abs{\norm{a_n+a}_E-2\norm{a}_E}\leq\norm{a_n-a}<\tau$ for each $n$. Consequently,
\begin{equation*}
\abs{\norm{y_n+y}_E-2\norm{y}_E}<3\tau \ \ \forall n\in \N.
\end{equation*}
Since $\norm{y_n/\norm{y_n}_E-y_n/\norm{y}_E}_E=\abs{1-\norm{y_n}_E/\norm{y}_E}<\tau/\norm{y}_E$ we get
\begin{equation}\label{eq:8}
\abs*{2-\norm*{\frac{y}{\norm{y}_E}+\frac{y_n}{\norm{y_n}_E}}_E}<\frac{4\tau}{\norm{y}_E}<\frac{2\tau}{\gamma} \ \ \forall n\in \N,
\end{equation}
where the last inequality holds because of \eqref{eq:6}. Note that $(x_{n,i})_{n\in \N}$ converges weakly to zero in $X_i$ for each $i\in I$ and
thus, by the representation of the dual of $\bigl[\bigoplus_{i\in J}X_i\bigr]_E$ as $\bigl[\bigoplus_{i\in J}X_i^*\bigr]_{E^{\prime}}$ and finiteness 
of $J$, the sequence $(y_n/\norm{y_n}_E)_{n\in \N}$ is also a weakly null sequence (as noted above, $(\norm{y_n})_{n\in \N}$ is bounded away from zero).\par
So from \eqref{eq:8} and the definition of $\alpha$ it follows that $\alpha\geq 2(1-\tau/\gamma)$. But $\gamma=\sqrt{\tau}$ and $\tau<(1-\alpha/2)^2$,
thus $2(1-\tau/\gamma)>\alpha$ and with this contradiction the proof is finished.
\end{Proof}

The above theorem reduces the case of infinite sums to the one of finite sums. The condition $\alpha<2$ is clearly necessary for 
$R(X)<2$. Unfortunately, the author does not know whether the simpler condition $\beta:=\sup_{i\in I}R(X_i)<2$ would be already enough to 
ensure that $\alpha<2$. The proof of \cite{dhompongsa}*{Theorem 7} shows that for $\beta<2$ one has for every finite subset $J\ssq I$ 
with $\abs{J}=N$ that $R(\bigl[\bigoplus_{i\in J}X_i\bigr]_E)\leq 2-\delta$, where first $\eps>0$ is chosen such that $\beta(1+N\eps)<2$ 
and then $0<\delta<\min\set*{2\delta_E(\eps),2-\beta(1+N\eps)}$, so it still might be that $R(\bigl[\bigoplus_{i\in J}X_i\bigr]_E)$ tends
to $2$ for $N\to \infty$.\par
Next we will discuss some applications of Theorem \ref{thm:garcia-falset}. First, since the Schur property is inherited by finite sums, 
we get the following corollary.
 
\begin{corollary}\label{cor:schur}
If $(X_i)_{i\in I}$ is a family of Banach spaces with the Schur property (in particular, a family of finite-dimensional Banach spaces)
and $\lin\set*{e_i:i\in I}$ is dense in $E$ with $\delta_E(1/4)>0$, then $R\paren*{\bigl[\bigoplus_{i\in I}X_i\bigr]_E}<2$. In particular,
$R\paren*{\bigl[\bigoplus_{i\in I}X_i\bigr]_p}<2$ for all $1<p<\infty$.
\end{corollary}

For another application of Theorem \ref{thm:garcia-falset} consider the following example.
\begin{example}\label{ex:c_0}
If $N\geq 2$ and $I_1, \dots, I_N$ are non-empty sets at least one of which is infinite, then
\begin{equation}\label{eq:9}
R\paren*{\Bigl[\bigoplus_{k=1}^{N}c_0(I_k)\Bigr]_p}=2^{1/p}
\end{equation}
for every $1\leq p<\infty$. Consequently, by Theorem \ref{thm:garcia-falset}, if $(I_k)_{k\in I}$ is any family of non-empty sets
we have that
\begin{equation*}
R\paren*{\Bigl[\bigoplus_{k\in I}c_0(I_k)\Bigr]_p}<2 \ \ \text{for} \ 1<p<\infty.
\end{equation*}
\end{example}

\begin{Proof}
To prove \eqref{eq:9} put $X:=\bigl[\bigoplus_{k=1}^{N}c_0(I_k)\bigr]_p$ and suppose without loss of generality that $I_1$ is infinite.
Fix a sequence $(i_n)_{n\in \N}$ of distinct elements of $I_1$ and any $j\in I_2$ and put $x_n:=(e_{i_n},0, \dots,0)\in S_X$ as well as 
$x:=(0, e_j, 0,\dots, 0)\in S_X$. Then $x_n\to 0$ weakly in $X$ and $\norm{x_n+x}_p=2^{1/p}$ for each $n$, thus $2^{1/p}\leq R(X)$.\par
To prove the reverse inequality let $x_n=(x_{n,1}, \dots, x_{n,N})\in B_X$ for each $n\in \N$ such that $x_n\to 0$ weakly and let 
$x=(x_1, \dots, x_N)\in B_X$. Without loss of generality we can suppose that $\lim_{n\to \infty}\norm{x_n+x}_p$ and also 
$a_k:=\lim_{n\to \infty}\norm{x_{n,k}}_{\infty}$ exists for each $k\in \set*{1, \dots, N}$.\par
Take an arbitrary $\eps>0$. Then for each $k\in \set*{1, \dots, N}$ the set $J_k:=\set*{i\in I_k:\abs{x_k(i)}>\eps}$ is finite. Since 
$x_n\to 0$ weakly we have $x_{n,k}(i)\to 0$ for all $k\in \set*{1, \dots, N}$ and all $i\in I_k$. It follows that there exists $n_0\in \N$ 
such that $\abs{x_{n,k}(i)}\leq\eps$ for all $k\in \set*{1, \dots, N}$, all $i\in J_k$ and all $n\geq n_0$.\par
But then $\abs{x_{n,k}(i)+x_k(i)}\leq\abs{x_{n,k}(i)}+\abs{x_k(i)}\leq\max\set*{\abs{x_{n,k}(i)},\abs{x_k(i)}}+\eps$ for all 
$k\in \set*{1, \dots, N}$, all $i\in I_k$ and all $n\geq n_0$. From this we can conclude
\begin{equation*}
\norm{x_n+x}_p^p=\sum_{k=1}^N\norm{x_{n,k}+x_k}_{\infty}^p\leq\sum_{k=1}^N(\max\set*{\norm{x_{n,k}}_{\infty},\norm{x_k}_{\infty}}+\eps)^p \ \ \forall n\geq n_0.
\end{equation*}
For $n\to \infty$ it follows that
\begin{equation*}
\lim_{n\to \infty}\norm{x_n+x}_p^p\leq\sum_{k=1}^N(\max\set*{a_k,\norm{x_k}_{\infty}}+\eps)^p.
\end{equation*}
Letting $\eps\to 0$ we obtain
\begin{align*}
&\lim_{n\to \infty}\norm{x_n+x}_p^p\leq\sum_{k=1}^N\max\set*{a_k^p,\norm{x_k}_{\infty}^p}\leq\sum_{k=1}^N(a_k^p+\norm{x_k}_{\infty}^p) \\
&=\lim_{n\to \infty}\norm{x_n}_p^p+\norm{x}_p^p\leq 2.
\end{align*}
Hence $\lim_{n\to \infty}\norm{x_n+x}_p\leq 2^{1/p}$ and we are done.
\end{Proof}

Recall that a Banach space $X$ is said to be a $U$-space if for any two sequences $(x_n)_{n\in \N}$ and $(y_n)_{n\in \N}$ in $S_X$ 
and every sequence $(x_n^*)_{n\in \N}$ in $S_{X^*}$ the conditions $\norm{x_n+y_n}\to 2$ and $x_n^*(x_n)=1$ for each $n\in \N$ imply 
$x_n^*(y_n)\to 1$. $U$-spaces were introduced by Lau in \cite{lau}. Uniformly rotund and uniformly smooth spaces are examples of 
$U$-spaces. Gao \cite{gao} defined the modulus of $u$-convexity of $X$ by 
\begin{equation*}
u_X(\eps):=\inf\set*{1-\norm*{\frac{x+y}{2}}:x,y\in S_X \ \exists x^*\in S_{X^*} \ x^*(x)=1, x^*(y)\leq 1-\eps}
\end{equation*}
for $0<\eps\leq 2$. Then $X$ is a $U$-space if and only if $u_X(\eps)>0$ for all $0<\eps\leq 2$. Obviously, $u_X\geq\delta_X$. 
By \cite{mazcunan-navvarro}*{Theorem 5} $R(X)<2$ if $u_X(\eps)>0$ for some $0<\eps<1$.\par
Putting several results together it is now possible to obtain the following corollary.
\begin{corollary}\label{cor:p-sum}
Let $(X_i)_{i\in I}$ be a family of Banach spaces and $1<p<\infty$. Suppose that there exist four pairwise disjoint (possibly empty) subsets 
$I_1, I_2, I_3, I_4\ssq I$ such that
\begin{enumerate}[\upshape(i)]
\item $X_i$ has the Schur property for each $i\in I_1$,
\item $\inf_{i\in I_2}u_{X_i}(\eps)>0$ for all $0<\eps\leq 2$,
\item for each $i\in I_3$ there is a set $J_i$ with $X_i=c_0(J_i)$.
\item $I_4$ is finite and $R(X_i)<2$ for all $i\in I_4$.
\end{enumerate}
Then $R(\bigl[\bigoplus_{i\in I}X_i\bigr]_p)<2$.
\end{corollary}

\begin{Proof}
Let us put $X:=\bigl[\bigoplus_{i\in I}X_i\bigr]_p$ and $X_k:=\bigl[\bigoplus_{i\in I_k}X_i\bigr]_p$ for $k=1,2,3,4$ (or $X_k=\set*{0}$ 
if $I_k=\emptyset$). By Corollary \ref{cor:schur} we have $R(X_1)<2$ and by Example \ref{ex:c_0} we have $R(X_3)<2$. Also, by 
\cite{hardtke}*{Corollary 3.17} and the remarks after \cite{hardtke}*{Definition 1.5} $X_2$ is again a $U$-space, so $R(X_2)<2$.
From the aforementioned result \cite{dhompongsa}*{Theorem 7} it follows that $R(X_4)<2$ and since $X\cong X_1\oplus_pX_2\oplus_pX_3\oplus_pX_4$, 
\cite{dhompongsa}*{Theorem 7} implies that $R(X)<2$.
\end{Proof}

The case of $c_0$-sums is not covered by the above results. However, it is easy to prove the following proposition directly.
\begin{proposition}\label{prop:c_0-sum}
Let $(X_i)_{i\in I}$ be any family of Banach spaces and $X:=\bigl[\bigoplus_{i\in I}X_i\bigr]_{c_0(I)}$. Then
\begin{equation*}
R(X)=\sup_{i\in I}R(X_i).
\end{equation*}
\end{proposition}

\begin{Proof}
We clearly have $\alpha:=\sup_{i\in I}R(X_i)\leq R(X)$. To prove the reverse inequality, fix any weakly null sequence $(x_n)_{n\in \N}=
((x_{n,i})_{i\in I})_{n\in \N}$ in $B_X$ and any $x=(x_i)_{i\in I}\in B_X$. Without loss of generality, we may assume that 
$\lim_{n\to \infty}\norm{x_n+x}_{\infty}$ exists.\par
Let $\eps>0$ be arbitrary. Then $J:=\set*{i\in I:\norm{x_i}\geq\eps}$ is finite, so by passing to an appropriate subsequence once more we may 
also assume that $\lim_{n\to \infty}\norm{x_{n,i}+x_i}$ exists for all $i\in J$.\par
Since $x_{n,i}\to 0$ weakly for all $i\in I$ it follows that $\lim_{n\to \infty}\norm{x_{n,i}+x_i}\leq R(X_i)\leq\alpha$ for all $i\in J$, so
$\norm{x_{n,i}+x_i}\leq\alpha+\eps$ for all $§i\in J$ and all sufficiently large $n$. But for $i\in I\sm J$ we have $\norm{x_{n,i}+x_i}\leq
\norm{x_{n,i}}+\norm{x_i}\leq1+\eps\leq\alpha+\eps$. Consequently, $\norm{x_n+x}_{\infty}\leq\alpha+\eps$ for all sufficiently large $n$, hence
$\lim_{n\to \infty}\norm{x_n+x}_{\infty}\leq\alpha+\eps$. Since $\eps>0$ was arbitrary, we are done.
\end{Proof}

Concerning $\ell^1$-sums it was already proved in \cite{kato}*{Theorem 3.13} that $R(X\oplus_1Y)<2$ if and only if both $X$ and $Y$ have the Schur property.
The proof of the ``only if'' part directly generalises to sums of arbitrarily many spaces and since it was proved in \cite{tanbay} that the $\ell^1$-sum
of any family of Banach spaces has the Schur property if and only if each summand has the Schur property, we obtain the following characterisation.
\begin{proposition}\label{prop:l_1}
Let $I$ be any index set with at least two elements. Let $(X_i)_{i\in I}$ be a family of Banach spaces and $X:=\bigl[\bigoplus_{i\in I}X_i\bigr]_1$. 
The following assertions are equivalent:
\begin{enumerate}[\upshape(i)]
\item $R(X)<2$,
\item $X_i$ has the Schur property for each $i\in I$,
\item $X$ has the Schur property,
\item $R(X)=1$.
\end{enumerate}
\end{proposition}

\subsection{Opial properties of finite absolute sums}\label{sub:finiteopial}
In this subsection we will briefly consider Opial properties of finite sums. This is surely well-known, 
but we will include the results and some of their proofs here as the author was not able to find them 
explicitly in the literature.\par
Recall that an absolute, normalised norm $\norm{\cdot}_E$ on $\R^m$ is said to be strictly monotone 
if for all $a=(a_1,\dots,a_m), b=(b_1,\dots,b_m)\in \R^m$ we have
\begin{equation*}
\norm{a}_E=\norm{b}_E \ \mathrm{and} \ \abs*{a_i}\leq\abs*{b_i} \ \forall i=1,\dots,m \ \Rightarrow \ \abs*{a_i}=\abs*{b_i} \ \forall i=1,\dots,m.
\end{equation*}
It is easy to see that strictly convex, absolute, normalised norms are strictly monotone.
\begin{proposition}\label{prop:finiteopial}
Let $\norm{\cdot}_E$ be an absolute, normalised norm on $\R^m$ and $X_1,\dots,X_m$ Banach spaces with 
nonstrict Opial property. Then $\bigl[\bigoplus_{i=1}^mX_i\bigr]_E$ has the nonstrict Opial property.
If moreover $\norm{\cdot}_E$ is strictly monotone and each $X_i$ has the Opial property, then 
$\bigl[\bigoplus_{i=1}^mX_i\bigr]_E$ also has the Opial property.
\end{proposition}
The proof is straightforward and will be omitted.\par
As is well-known, every strictly monotone, absolute, normalised norm on $\R^m$ is actually uniformly monotone 
in the following sense (the proof consists in an easy compactness argument).
\begin{lemma}\label{lemma:unifmonotone}
Let $\norm{\cdot}_E$ be a strictly monotone, absolute, normalised norm on $\R^m$. Let $\eps,R>0$.
The there exists $\delta>0$ such that for all $a=(a_1,\dots,a_m), b=(b_1,\dots,b_m)\in \R^m$ with
$\norm{b}_E\leq R$ and $\abs*{a_i}\leq\abs*{b_i}$ for $i=1\dots,m$ we have
\begin{equation*}
\norm{b}_E-\norm{a}_E<\delta \ \Rightarrow \ \abs*{b_i}-\abs*{a_i}<\eps \ \forall i=1,\dots,m.
\end{equation*}
\end{lemma}
Utilizing this fact, one can see the following.
\begin{proposition}\label{prop:finiteunifopial}
Let $\norm{\cdot}_E$ be an absolute, normalised norm on $\R^m$ which is strictly monotone and $X_1,\dots,X_m$ Banach 
spaces with the uniform Opial property. Then $X:=\bigl[\bigoplus_{i=1}^mX_i\bigr]_E$ also has the uniform Opial property.
\end{proposition}

\begin{Proof}
Let $\eps,R>0$ and put $\eta:=\min\set*{\eta_{X_i}(\eps/m,R):i=1,\dots,m}$. Choose a $0<\delta\leq1$ according to Lemma 
\ref{lemma:unifmonotone} corresponding to the values $\eta$ and $3R+1$.\par
Now consider a weakly nullsequence $(x_n)_{n\in \N}=((x_{n,1},\dots,x_{n,m}))_{n\in \N}$ in $X$ with $\limsup\norm{x_n}_E\leq R$ and an element 
$y=(y_1,\dots,y_m)\in X$ with $\norm{y}_E\geq\eps$. Since $\norm{y}_E\leq\sum_{i=1}^m\norm{y_i}$ there is some $i_0\in\set*{1,\dots,m}$ with
$\norm{y_{i_0}}\geq\eps/m$. There is no loss of generality in assuming that all the limits in the following calculations exist. From the 
definition of $\eta$ we get
\begin{equation*}
\lim_{n\to \infty}\norm{x_{n,i_0}}+\eta\leq\lim_{n\to \infty}\norm{x_{n,i_0}-y_{i_0}}.
\end{equation*}
Since each $X_i$ has in particular the nonstrict Opial property, we also have
\begin{equation*}
\lim_{n\to \infty}\norm{x_{n,i}}\leq\lim_{n\to \infty}\norm{x_{n,i}-y_i} \ \ \forall i\in\set*{1,\dots,m}\sm\set*{i_0}.
\end{equation*}
If $\norm{y}_E\leq 2R+1$, then $\lim_{n\to \infty}\norm{x_n-y}_E\leq\lim_{n\to \infty}\norm{x_n}_E+2R+1\leq3R+1$ and the choice 
of $\delta$ implies $\lim_{n\to \infty}\norm{x_n}_E+\delta\leq\lim_{n\to \infty}\norm{x_n-y}_E$.\par
If on the other hand $\norm{y}_E>2R+1$, then $\lim_{n\to \infty}\norm{x_n-y}_E\geq\norm{y}_E-\lim_{n\to \infty}\norm{x_n}_E\geq R+1\geq
\lim_{n\to \infty}\norm{x_n}_E+\delta$. So $X$ has the uniform Opial property.
\end{Proof}

\subsection{Opial properties of some infinite sums}\label{sub:infiniteopial}
We will first show that the Opial and nonstrict Opial property are preserved under infinite $\ell^p$-sums.
\begin{proposition} \label{prop:p-opial}
If $1\leq p<\infty$, $I$ is any index set and $(X_i)_{i\in I}$ a family of Banach spaces with the Opial property (nonstrict Opial property),
then $X:=\bigl[\bigoplus_{i\in I}X_i\bigr]_p$ also has the Opial property (nonstrict Opial property).
\end{proposition}

\begin{Proof}
We will only prove the strict case, the nonstrict case is treated analogously. Let $x_n=(x_{n,i})_{i\in I}\in X$ for every $n\in \N$ 
such that $(x_n)_{n\in \N}$ converges weakly to zero and let $x=(x_i)_{i\in I}\in X\sm\set*{0}$. Fix $i_0\in I$ with $x_{i_0}\neq 0$.
We may assume that $\lim_{n\to \infty}\norm{x_n}_p$ and $\lim_{n\to \infty}\norm{x_n-x}_p$ as well as $a:=\lim_{n\to \infty}\norm{x_{n,i_0}}$
and $b:=\lim_{n\to \infty}\norm{x_{n,i_0}-x_{i_0}}$ exist. Note also that $(x_{n,i})_{n\in \N}$ is a weakly nullsequence in $X_i$ for 
each $i\in I$. So since $X_{i_0}$ has the Opial property it follows that $\delta:=b^p-a^p>0$. Put $K:=\sup_{n\in \N}\norm{x_n}_p$ and
let $0<\eps\leq1$. We can find a finite set $J\ssq I$ with $i_0\in J$ such that
\begin{equation}\label{eq:10}
\norm*{(\norm{x_i}\chi_{I\sm J}(i))_{i\in I}}_p\leq\eps,
\end{equation}
where $\chi_{I\sm J}$ denotes the characteristic function of $I\sm J$. By passing to a further subsequence, we can assume that
$\lim_{n\to \infty}\norm{x_{n,i}}$ and $\lim_{n\to \infty}\norm{x_{n,i}-x_i}$ exist for all $i\in J$. Then, using the Opial property of 
each of the summands $X_i$, the definition of $\delta$ and \eqref{eq:10}, we obtain
\begin{align*}
&\lim_{n\to \infty}\norm{x_n}_p^p=\sum_{i\in J\sm\set*{i_0}}\lim_{n\to \infty}\norm{x_{n,i}}^p+a^p+\lim_{n\to \infty}
\norm*{(\norm{x_{n,i}}\chi_{I\sm J}(i))_{i\in I}}_p^p \\
&\leq\lim_{n\to \infty}\sum_{i\in J\sm\set*{i_0}}\norm{x_{n,i}-x_i}^p+b^p-\delta+\lim_{n\to \infty}\norm*{(\norm{x_{n,i}}\chi_{I\sm J}(i))_{i\in I}}_p^p \\
&\leq\lim_{n\to \infty}\sum_{i\in J}\norm{x_{n,i}-x_i}^p-\delta+\lim_{n\to \infty}\paren*{\norm*{(\norm{x_{n,i}-x_i}\chi_{I\sm J}(i))_{i\in I}}_p+\eps}^p.
\end{align*}
But, since $\abs*{s^p-t^p}\leq pA^{p-1}\abs*{s-t}$ for all $0\leq s,t\leq A$, we also have
\begin{align*}
&\lim_{n\to \infty}\paren*{\norm*{(\norm{x_{n,i}-x_i}\chi_{I\sm J}(i))_{i\in I}}_p+\eps}^p\leq\lim_{n\to \infty}\norm*{(\norm{x_{n,i}-x_i}\chi_{I\sm J}(i))_{i\in I}}_p^p \\
&+\lim_{n\to \infty}\abs*{\paren*{\norm*{(\norm{x_{n,i}-x_i}\chi_{I\sm J}(i))_{i\in I}}_p+\eps}^p-\norm*{(\norm{x_{n,i}-x_i}\chi_{I\sm J}(i))_{i\in I}}_p^p} \\
&\leq\lim_{n\to \infty}\norm*{(\norm{x_{n,i}-x_i}\chi_{I\sm J}(i))_{i\in I}}_p^p+p(K+\norm{x}_p+1)^{p-1}\eps.
\end{align*}
It follows that
\begin{equation*}
\lim_{n\to \infty}\norm{x_n}_p^p\leq\lim_{n\to \infty}\norm{x_n-x}_p^p-\delta+p(K+\norm{x}_p+1)^{p-1}\eps.
\end{equation*}
Since $\eps\in (0,1]$ was arbitrary and $\delta$ independent of $\eps$, we conclude
\begin{equation*}
\lim_{n\to \infty}\norm{x_n}_p^p\leq\lim_{n\to \infty}\norm{x_n-x}_p^p-\delta<\lim_{n\to \infty}\norm{x_n-x}_p^p
\end{equation*}
and the proof is finished.
\end{Proof}

$c_0$ is a typical example of a Banach space which has the nonstrict Opial property but not the usual (strict) Opial property.
Next we will see that $c_0$-sums preserve the nonstrict Opial property.
\begin{proposition}\label{prop:opial-c_0}
Let $I$ be any index set and $(X_i)_{i\in I}$ a family of Banach spaces with the nonstrict Opial property.
Then $X:=\bigl[\bigoplus_{i\in I}X_i\bigr]_{c_0(I)}$ has the nonstrict Opial property.
\end{proposition}

\begin{Proof}
Let $x_n=(x_{n,i})_{i\in I}\in X$ for every $n\in \N$ such that $(x_n)_{n\in \N}$ converges weakly to zero and let 
$x=(x_i)_{i\in I}\in X$. Take $\eps>0$ to be arbitrary and find a finite subset $J\ssq I$ such that $\norm{x_i}\leq\eps$
for every $i\in I\sm J$. Again there is no loss of generality in assuming that all the limits involved in the following 
calculations exist. Since each $X_i$ has the nonstrict Opial property, we have
\begin{equation*}
\lim_{n\to \infty}\norm{x_{n,i}}\leq\lim_{n\to \infty}\norm{x_{n,i}-x_i} \ \ \forall i\in J.
\end{equation*}
Therefore we obtain
\begin{align*}
&\lim_{n\to \infty}\norm{x_n}_{\infty}=\max\set*{\max_{i\in J}\lim_{n\to \infty}\norm{x_{n,i}},\lim_{n\to \infty}\norm{(\norm{x_{n,i}}\chi_{I\sm J}(i))_{i\in I}}_{\infty}} \\
&\leq\max\set*{\max_{i\in J}\lim_{n\to \infty}\norm{x_{n,i}-x_i},\lim_{n\to \infty}\norm{(\norm{x_{n,i}}\chi_{I\sm J}(i))_{i\in I}}_{\infty}} \\
&\leq\lim_{n\to \infty}\max\set*{\max_{i\in J}\norm{x_{n,i}-x_i},\norm{(\norm{x_{n,i}-x_i}\chi_{I\sm J}(i))_{i\in I}}_{\infty}}+\eps \\
&=\lim_{n\to \infty}\norm{x_n-x}_{\infty}+\eps
\end{align*}
and since $\eps>0$ was arbitrary we are done.
\end{Proof}

Concerning the uniform Opial property, we have the following result for infinite $\ell^p$-sums, resembling in structure Theorem \ref{thm:garcia-falset}.
\begin{theorem}\label{thm:unifopial}
Let $1\leq p<\infty$ and let $I$ be an infinite index set. For a family $(X_i)_{i\in I}$ of Banach spaces put $X_J:=\bigl[\bigoplus_{i\in J}X_i\bigr]_p$
for every finite $J\ssq I$. Suppose that
\begin{equation*}
\omega(\eps,R):=\inf\set*{\eta_{X_J}(\eps,R):J\ssq I \ \mathrm{finite}}>0 \ \ \forall \eps,R>0.
\end{equation*}
Then $X:=\bigl[\bigoplus_{i\in I}X_i\bigr]_p$ has the uniform Opial property.
\end{theorem}

\begin{Proof}
Let $0<\eps\leq 1$ and $R>0$. We put $\nu:=\min\set*{3R+1,\omega(\eps/2,R)}$ and $\tau:=\min\set*{1,3R+1-((3R+1)^p-\nu^p)^{1/p}}$.
Now let us consider a weakly nullsequence $(x_n)_{n\in \N}=((x_{n,i})_{i\in I})_{n\in \N}$ in $X$ with $\limsup\norm{x_n}_p\leq R$ 
and let $x=(x_i)_{i\in I}\in X$ with $\norm{x}_p\geq\eps$. As before, we may assume that $\lim_{n\to \infty}\norm{x_n}_p$ and 
$\lim_{n\to \infty}\norm{x_n-x}_p$ exist. Let $K:=\sup_{n\in \N}\norm{x_n}_p$. For $0<\alpha\leq\eps/2$ we can find a finite 
subset $J\ssq I$ such that
\begin{equation*}
\norm*{(\norm{x_i}\chi_{I\sm J}(i))_{i\in I}}_p\leq\alpha.
\end{equation*}
It follows that
\begin{equation}\label{eq:11}
\norm*{\sum_{i\in J}\norm{x_i}e_i}_p\geq\norm{x}_p-\alpha\geq\eps/2.
\end{equation}
We may also assume that $\lim_{n\to \infty}\norm{x_{n,i}}$ and $\lim_{n\to \infty}\norm{x_{n,i}-x_i}$ exist for all $i\in J$.
Analogous to the proof of Proposition \ref{prop:p-opial} we can show that
\begin{align}\label{eq:12}
&\lim_{n\to \infty}\norm{x_n}_p^p\leq\lim_{n\to \infty}\sum_{i\in J}\norm{x_{n,i}}^p+\lim_{n\to \infty}\norm*{(\norm{x_{n,i}-x_i}\chi_{I\sm J}(i))_{i\in I}}_p^p \nonumber \\
&+p(K+\norm{x}_p+1)^{p-1}\alpha.
\end{align}
If we put $y_n:=(x_{n,i})_{i\in J}$ for each $n\in \N$ and $y:=(x_i)_{i\in J}$, then $(y_n)_{n\in \N}$ is a weakly nullsequence in $X_J$ 
with $\lim_{n\to \infty}\norm{y_n}_p\leq\lim_{n\to \infty}\norm{x_n}_p\leq R$ and $y\in X_J$ with $\norm{y}_p\geq\eps/2$ (because of \eqref{eq:11}), thus
\begin{equation}\label{eq:13}
\lim_{n\to \infty}\norm{y_n}_p+\eta_{X_J}(\eps/2,R)\leq\lim_{n\to \infty}\norm{y_n-y}_p.
\end{equation}
Since $(a-b)^p\leq a^p-b^p$ for all $a\geq b\geq 0$ we can deduce from \eqref{eq:12} and \eqref{eq:13} that
\begin{align*}
&\lim_{n\to \infty}\norm{x_n}_p^p\leq\lim_{n\to \infty}\norm{y_n-y}_p^p-\eta_{X_J}(\eps/2,R)^p \\
&+\lim_{n\to \infty}\norm*{(\norm{x_{n,i}-x_i}\chi_{I\sm J}(i))_{i\in I}}_p^p+p(K+\norm{x}_p+1)^{p-1}\alpha \\
&\leq\lim_{n\to \infty}\norm{x_n-x}_p^p-\nu^p+p(K+\norm{x}_p+1)^{p-1}\alpha.
\end{align*}
Letting $\alpha\to 0$ we obtain
\begin{equation}\label{eq:14}
\lim_{n\to \infty}\norm{x_n}_p^p\leq\lim_{n\to \infty}\norm{x_n-x}_p^p-\nu^p.
\end{equation}
If $\norm{x}_p\geq 2R+1$ then $\lim_{n\to \infty}\norm{x_n-x}_p\geq2R+1-\lim_{n\to \infty}\norm{x_n}_p\geq
\lim_{n\to \infty}\norm{x_n}_p+1\geq\lim_{n\to \infty}\norm{x_n}_p+\tau$.\par
Now consider the case $\norm{x}_p<2R+1$. Define $f(s):=s-(s^p-\nu^p)^{1/p}$ for all $s\geq\nu$. It is easily checked that $f$ is 
decreasing. Since $\lim_{n\to \infty}\norm{x_n-x}_p\leq \lim_{n\to \infty}+\norm{x}_p\leq3R+1$ it follows from \eqref{eq:14} that
\begin{align*}
&\lim_{n\to \infty}\norm{x_n}_p\leq\lim_{n\to \infty}\norm{x_n-x}_p-f(\lim_{n\to \infty}\norm{x_n-x}_p) \\
&\leq\lim_{n\to \infty}\norm{x_n-x}_p-f(3R+1)\leq\lim_{n\to \infty}\norm{x_n-x}_p-\tau
\end{align*}
and the proof is complete.
\end{Proof}

As a corollary we obtain again the already known result that the $\ell^p$-sum of any family of Banach spaces with the Schur property
has the uniform Opial property (see \cite{prus}*{Example 4.23 (2.)} or \cite{saejung}*{Theorem 7}).
\begin{corollary}\label{cor:opialschur}
Let $1\leq p<\infty$ and let $(X_i)_{i\in I}$ be a family of Banach spaces with the Schur property. Then $\bigl[\bigoplus_{i\in I}X_i\bigr]_p$ 
has the uniform Opial property.
\end{corollary}

The author does not know whether the conditon $\inf_{i\in I}\eta_{X_i}>0$ is already enough to ensure that $\bigl[\bigoplus_{i\in I}X_i\bigr]_p$
has the uniform Opial property (the proof of Proposition \ref{prop:finiteunifopial} does not give a uniform lower bound for the moduli of the 
finite sums).

\section{Some Opial-type properties in Lebesgue-Bochner spaces}
We consider a complete, finite measure space $(S,\mathcal{A},\mu)$ and a real Banach space $X$. First recall that for 
$1\leq p\leq\infty$ the Lebesgue-Bochner space $L^p(\mu,X)$ is defined as the space of all Bochner-measurable functions 
$f:S \rightarrow X$ (modulu equality $\mu$-almost everywhere) such that $\norm{f(\cdot)}\in L^p(\mu)$. Equipped with 
the norm $\norm{f}_p:=\norm{\norm{f(\cdot)}}_p$, $L^p(\mu,X)$ becomes a Banach space.\par
As was mentioned in the introduction, even the spaces $L^p[0,1]$, $1<p<\infty, p\neq 2$, of scalar-valued functions do not 
have the Opial property. However, some results which are in a certain sense analogous to the Opial property are available.
For example it was shown in \cite{brezis} that any bounded sequence $(f_n)_{n\in \N}$ in $L^p(\mu)$ ($0<p<\infty$) which 
converges pointwise almost everywhere to a function $f\in L^p(\mu)$ satisfies
\begin{equation*}
\lim_{n\to \infty}\paren*{\norm{f_n}_p^p-\norm{f_n-f}_p^p}=\norm{f}_p^p
\end{equation*}
and hence
\begin{equation*}
\liminf_{n\to \infty}\norm{g_n-g}_p^p=\liminf_{n\to \infty}\norm{g_n}_p^p+\norm{g}_p^p
\end{equation*}
for any bounded sequence $(g_n)_{n\in \N}$ in $L^p(\mu)$ which converges pointwise almost everywhere to zero 
and every $g\in L^p(\mu)$.\par
In \cite{besbes}*{Chapter 2, Lemma 3.3} it was shown that any sequence $(f_n)_{n\in \N}$ in $L^1(\mu,X)$ (where 
$(S,\mathcal{A},\mu)$ is a probability space and $X$ an arbitrary Banach space) and any $f\in L^1(\mu,X)$ such that
\begin{equation*}
\lim_{n\to \infty}\mu(\set*{t\in S:\norm{f_n(t)-f(t)}\geq\eps})=0 \ \ \forall\eps>0
\end{equation*}
satisfy the equality
\begin{equation*}
\liminf_{n\to \infty}\norm{f_n-f}_1+\norm{f-g}_1=\liminf_{n\to \infty}\norm{f_n-g}_1 
\end{equation*}
for every $g\in L^1(\mu,X)$.\par
We are now going to consider pointwise weak convergence almost everywhere in Lebesgue-Bochner spaces and prove some results 
analogous to the Opial property in this setting.
\begin{proposition}\label{prop:opialBochnerLp}
Let $(S,\mathcal{A},\mu)$ be a complete, finite measure space, $1\leq p<\infty$ and $X$ a Banach space with the nonstrict
Opial property. Let $(f_n)_{n\in \N}$ be a bounded sequence in $L^p(\mu,X)$ such that $(f_n(t))_{n\in \N}$ converges weakly
to zero for almost every $t\in S$. Suppose further that there is a function $g\in L^p(\mu)$ such that $\norm{f_n(t)}\to g(t)$
for almost every $t\in S$. Then
\begin{align*}
&\int_A\liminf_{n\to \infty}\norm{f_n(t)-f(t)}^p\,\mathrm{d}\mu(t)-\int_Ag(t)^p\,\mathrm{d}\mu(t) \\
&\leq\limsup_{n\to \infty}\norm{f_n-f}_p^p-\limsup_{n\to \infty}\norm{f_n}_p^p
\end{align*}
holds for every $f\in L^p(\mu,X)$ and every $A\in \mathcal{A}$. In particular,
\begin{equation*}
\limsup_{n\to \infty}\norm{f_n}_p\leq\limsup_{n\to \infty}\norm{f_n-f}_p \ \ \forall f\in L^p(\mu,X).
\end{equation*}
\end{proposition}

\begin{Proof}
Without loss of generality we can assume that $\lim_{n\to \infty}\norm{f_n(t)}=g(t)$ and $f_n(t)\to 0$ weakly 
for every $t\in S$ and also that $\lim_{n\to \infty}\norm{f_n}_p$ and $\lim_{n\to \infty}\norm{f_n-f}_p$ 
exist. Now let $A\in \mathcal{A}, f\in L^p(\mu,X)$ and $0<\eps<1$. By the equi-integrability of finite subsets 
of $L^1(\mu)$ there exists $\delta>0$ such that 
\begin{align}\label{eq:15}
&B\in \mathcal{A},\,\mu(B)\leq\delta \ \Rightarrow \ \int_Bh(t)\,\mathrm{d}\mu(t)\leq\eps \\
&\mathrm{for\ each}\ h\in\set*{\norm{f(\cdot)}^p,g^p,\liminf_{n\to \infty}\norm{f_n(\cdot)-f(\cdot)}^p}. \nonumber
\end{align}
By Egorov's theorem there exists $C\in \mathcal{A}$ with $\mu(S\sm C)\leq\delta$ such that $\lim_{n\to \infty}\norm{f_n(t)}^p=g(t)^p$
uniformly in $t\in C$, which implies
\begin{equation}\label{eq:16}
\lim_{n\to \infty}\int_C\norm{f_n(t)}^p\,\mathrm{d}\mu(t)=\int_Cg(t)^p\,\mathrm{d}\mu(t).
\end{equation}
By passing to a further subsequence if necessary, we may assume that the limit
\begin{equation*}
\lim_{n\to \infty}\int_C\norm{f_n(t)-f(t)}^p\,\mathrm{d}\mu(t)
\end{equation*}
exists. Now we can calculate, using \eqref{eq:16} and the nonstrict Opial property of $X$,
\begin{align}\label{eq:17}
&\lim_{n\to \infty}\norm{f_n}_p^p=\int_Cg(t)^p\,\mathrm{d}\mu(t)+\lim_{n\to \infty}\int_{S\sm C}\norm{f_n(t)}^p\,\mathrm{d}\mu(t) \nonumber \\
&=\int_{C\cap A}g(t)^p\,\mathrm{d}\mu(t)+\int_{C\sm A}g(t)^p\,\mathrm{d}\mu(t)+\lim_{n\to \infty}\int_{S\sm C}\norm{f_n(t)}^p\,\mathrm{d}\mu(t) \nonumber \\
&\leq\int_{C\cap A}g(t)^p\,\mathrm{d}\mu(t)+\int_{C\sm A}\liminf_{n\to \infty}\norm{f_n(t)-f(t)}^p\,\mathrm{d}\mu(t) \nonumber \\
&+\lim_{n\to \infty}\int_{S\sm C}\norm{f_n(t)}^p\,\mathrm{d}\mu(t)\leq\int_{C}\liminf_{n\to \infty}\norm{f_n(t)-f(t)}^p\,\mathrm{d}\mu(t) \nonumber \\
&+\int_{C\cap A}\paren*{g(t)^p-\liminf_{n\to \infty}\norm{f_n(t)-f(t)}^p}\,\mathrm{d}\mu(t)+\lim_{n\to \infty}\int_{S\sm C}\norm{f_n(t)}^p\,\mathrm{d}\mu(t).
\end{align}
But
\begin{equation}\label{eq:18}
\abs*{\int_{C\cap A}g(t)^p\,\mathrm{d}\mu(t)-\int_{A}g(t)^p\,\mathrm{d}\mu(t)}=\int_{A\sm C}g(t)^p\,\mathrm{d}\mu(t)\leq\eps
\end{equation}
because of $\mu(S\sm C)\leq\delta$ and \eqref{eq:15}. Analogously,
\begin{equation}\label{eq:19}
\int_{A\sm C}\liminf_{n\to \infty}\norm{f_n(t)-f(t)}^p\,\mathrm{d}\mu(t)\leq\eps.
\end{equation}
Putting \eqref{eq:17}, \eqref{eq:18} and \eqref{eq:19} together we obtain
\begin{align}\label{eq:20}
&\lim_{n\to \infty}\norm{f_n}_p^p\leq\int_{C}\liminf_{n\to \infty}\norm{f_n(t)-f(t)}^p\,\mathrm{d}\mu(t)+\lim_{n\to \infty}\int_{S\sm C}\norm{f_n(t)}^p\,\mathrm{d}\mu(t) \nonumber \\
&+\int_{A}\paren*{g(t)^p-\liminf_{n\to \infty}\norm{f_n(t)-f(t)}^p}\,\mathrm{d}\mu(t)+2\eps \nonumber \\
&\leq\lim_{n\to \infty}\int_{C}\norm{f_n(t)-f(t)}^p\,\mathrm{d}\mu(t)+\lim_{n\to \infty}\int_{S\sm C}\norm{f_n(t)}^p\,\mathrm{d}\mu(t) \nonumber \\
&+\int_{A}\paren*{g(t)^p-\liminf_{n\to \infty}\norm{f_n(t)-f(t)}^p}\,\mathrm{d}\mu(t)+2\eps,
\end{align}
where we have used Fatou's lemma in the second step.\par
Since $\mu(S\sm C)\leq\delta$ we have $\int_{S\sm C}\norm{f(t)}^p\,\mathrm{d}\mu(t)\leq\eps$, by \eqref{eq:15}. Hence
\begin{equation*}
\lim_{n\to \infty}\int_{S\sm C}\norm{f_n(t)}^p\,\mathrm{d}\mu(t)\leq\lim_{n\to \infty}\paren*{\paren*{\int_{S\sm C}\norm{f_n(t)-f(t)}^p\,\mathrm{d}\mu(t)}^{1/p}+\eps^{1/p}}^p
\end{equation*}
Since $\abs*{s^p-t^p}\leq pA^{p-1}\abs*{s-t}$ for all $0\leq s,t\leq A$, we obtain as in the proof of Proposition \ref{prop:p-opial}
\begin{equation}\label{eq:21}
\lim_{n\to \infty}\int_{S\sm C}\norm{f_n(t)}^p\,\mathrm{d}\mu(t)\leq\lim_{n\to \infty}\int_{S\sm C}\norm{f_n(t)-f(t)}^p\,\mathrm{d}\mu(t)+pL^{p-1}\eps^{1/p},
\end{equation}
where $L:=\norm{f}_p+1+\sup_{n\in \N}\norm{f_n}_p$. From \eqref{eq:20} and \eqref{eq:21} it follows that
\begin{align*}
&\lim_{n\to \infty}\norm{f_n}_p^p\leq\lim_{n\to \infty}\int_{S}\norm{f_n(t)-f(t)}^p\,\mathrm{d}\mu(t)+pL^{p-1}\eps^{1/p} \\
&+\int_{A}\paren*{g(t)^p-\liminf_{n\to \infty}\norm{f_n(t)-f(t)}^p}\,\mathrm{d}\mu(t)+2\eps.
\end{align*}
Letting $\eps\to 0$ now leads to the desired inequality.
\end{Proof}

If $X$ has the Opial property, we have the following corollary.
\begin{corollary}\label{cor:opialBochnerLp}
Let $(S,\mathcal{A},\mu)$ be a complete, finite measure space, $1\leq p<\infty$ and $X$ a Banach space with the 
Opial property. Let $(f_n)_{n\in \N}$ be a bounded sequence in $L^p(\mu,X)$ such that $(f_n(t))_{n\in \N}$ converges 
weakly to zero for almost every $t\in S$. Suppose further that there is a function $g\in L^p(\mu)$ such that 
$\norm{f_n(t)}\to g(t)$ for almost every $t\in S$. Then
\begin{equation*}
\limsup_{n\to \infty}\norm{f_n}_p<\limsup_{n\to \infty}\norm{f_n-f}_p \ \ \forall f\in L^p(\mu,X)\sm\set*{0}.
\end{equation*}
\end{corollary}

\begin{Proof}
Just put $A:=\set*{t\in S:f(t)\neq 0}$ in Proposition \ref{prop:opialBochnerLp}. Then $\mu(A)>0$ and since $X$ has
the Opial property we have $\liminf\norm{f_n(t)-f(t)}<g(t)$ for every $t\in A$, so the result follows from Proposition 
\ref{prop:opialBochnerLp}.
\end{Proof}

In the case that $X$ even has the uniform Opial property, we have the following two results.
\begin{theorem}\label{thm:unifopialBochnerLp1}
Let $(S,\mathcal{A},\mu)$ be a complete, finite measure space, $1\leq p<\infty$ and $X$ a Banach space with the uniform
Opial property. Let $M,R>0$ and $f\in L^p(\mu,X)\sm\set*{0}$. Then there exists $\eta>0$ such that the following holds: 
whenever $(f_n)_{n\in \N}$ is a sequence in $L^p(\mu,X)$ with $\sup_{n\in \N}\norm{f_n}_p\leq R$ such that $(f_n(t))_{n\in \N}$ 
converges weakly to zero and $\lim_{n\to \infty}\norm{f_n(t)}\leq M$ for almost every $t\in S$, then
\begin{equation*}
\limsup_{n\to \infty}\norm{f_n}_p+\eta\leq\limsup_{n\to \infty}\norm{f_n-f}_p.
\end{equation*}
\end{theorem}

\begin{Proof}
We define $\tau:=\norm{f}_p(2\mu(S))^{-1/p}$ and $A:=\set*{t\in S:\norm{f(t)}\geq\tau}$. If $\mu(A)=0$, then we would obtain 
$\norm{f}_p^p\leq\mu(S\sm A)\tau^p\leq\norm{f}_p^p/2$, contradicting the fact that $f\in L^p(\mu,X)\sm\set*{0}$. Thus $\mu(A)>0$.\par
Next we put $w:=\eta_X(\tau,M)$, $\delta:=\min\set*{(3R+1)^p,\mu(A)w^p}$, $\omega:=§R+1-((§R+1)^p-\delta)^{1/p}$ and finally 
$\eta:=\min\set*{\omega,1}$.\par
Now let $(f_n)_{n\in \N}$ be as above. Without loss of generality we may assume that $g(t):=\lim_{n\to \infty}\norm{f_n(t)}\leq M$ 
and $f_n(t)\to 0$ weakly for every $t\in S$. The definition of $\eta_X$ implies
\begin{equation*}
\liminf_{n\to \infty}\norm{f_n(t)-f(t)}-g(t)\geq\eta_X(\tau,M)=w \ \ \forall t\in A.
\end{equation*}
Since $(a-b)^p\leq a^p-b^p$ for all $a\geq b\geq 0$, it follows that
\begin{equation*}
\liminf_{n\to \infty}\norm{f_n(t)-f(t)}^p-g(t)^p\geq w^p \ \ \forall t\in A.
\end{equation*}
Combinig this with Proposition \ref{prop:opialBochnerLp} leads to
\begin{equation*}
\limsup_{n\to \infty}\norm{f_n-f}_p^p-\limsup_{n\to \infty}\norm{f_n}_p^p\geq\mu(A)w^p\geq\delta.
\end{equation*}
As in the proof of Theorem \ref{thm:unifopial}, by distinguishing the two cases $\norm{f}_p\geq 2R+1$ and $\norm{f}_p<2R+1$, 
we can deduce from this that
\begin{equation*}
\limsup_{n\to \infty}\norm{f_n}_p+\eta\leq\limsup_{n\to \infty}\norm{f_n-f}_p.
\end{equation*}
\end{Proof}

\begin{theorem}\label{thm:unifopialBochnerLp2}
Let $(S,\mathcal{A},\mu)$ be a complete, finite measure space, $1\leq p<\infty$ and $X$ a Banach space with the uniform 
Opial property. Let $p<r\leq\infty$ and $\eps,M,R,K>0$. Then there exists $\eta>0$ such that the following holds: whenever
$(f_n)_{n\in \N}$ is a sequence in $L^p(\mu,X)$ with $\sup_{n\in \N}\norm{f_n}_p\leq R$ such that $(f_n(t))_{n\in \N}$
converges weakly to zero and $\lim_{n\to \infty}\norm{f_n(t)}\leq M$ for almost every $t\in S$ and $f\in L^r(\mu,X)\ssq L^p(\mu,X)$
such that $\norm{f}_r\leq K$ and $\norm{f}_p\geq\eps$, then
\begin{equation*}
\limsup_{n\to \infty}\norm{f_n}_p+\eta\leq\limsup_{n\to \infty}\norm{f_n-f}_p.
\end{equation*}
\end{theorem}

\begin{Proof}
We put $s:=r/p\in(1,\infty]$. Let $s^{\prime}\in [1,\infty)$ such that $1/s^{\prime}+1/s=1$. Choose $0<\tau<\eps\mu(S)^{-1/p}$ and 
put $Q:=(\eps^p-\mu(S)\tau^p)^{s^{\prime}}K^{-ps^{\prime}}$. Let $w:=\eta_X(\tau,M)$ and $\delta:=\min\set*{Qw^p,(3R+1)^p}$. $\omega$ 
and $\eta$ are also defined as in the previous proof.\par
Now let $(f_n)_{n\in \N}$ and $f$ be as above. For $A:=\set*{t\in S:\norm{f(t)}\geq\tau}$ we have
\begin{align*}
&\eps^p\leq\norm{f}_p^p=\int_A\norm{f(t)}^p\,\mathrm{d}\mu(t)+\int_{S\sm A}\norm{f(t)}^p\,\mathrm{d}\mu(t) \\
&\leq\int_A\norm{f(t)}^p\,\mathrm{d}\mu(t)+\mu(S\sm A)\tau^p\leq\mu(A)^{1/s^{\prime}}\norm{f}_r^{p}+\mu(S)\tau^p \\
&\leq\mu(A)^{1/s^{\prime}}K^p+\mu(S)\tau^p,
\end{align*}
where we have used H\"older's inequality in the second line. It follows that $\mu(A)\geq Q$. As in the previous proof we can deduce that
\begin{equation*}
\limsup_{n\to \infty}\norm{f_n-f}_p^p-\limsup_{n\to \infty}\norm{f_n}_p^p\geq\mu(A)w^p\geq Qw^p\geq\delta
\end{equation*}
and from there we get to 
\begin{equation*}
\limsup_{n\to \infty}\norm{f_n}_p+\eta\leq\limsup_{n\to \infty}\norm{f_n-f}_p
\end{equation*}
as in the proof of Theorem \ref{thm:unifopial}.
\end{Proof}

\address
\email

\end{document}